\documentclass[11pt]{article} 
\usepackage{amsmath}
\usepackage{amsfonts}
\usepackage[latin1]{inputenc}
\usepackage[T1]{fontenc}
\usepackage[french]{babel}

\makeatletter
\def\@sect#1#2#3#4#5#6[#7]#8{%
  \ifnum #2>\c@secnumdepth
    \let\@svsec\@empty
  \else
    \refstepcounter{#1}%
    \protected@edef\@svsec{\@seccntformat{#1}\relax}%
  \fi
  \@tempskipa #5\relax
  \ifdim \@tempskipa>\z@
    \begingroup
      #6{%
        \@hangfrom{\hskip #3\relax\@svsec}%
          \interlinepenalty \@M #8\@@par}%
    \endgroup
    \csname #1mark\endcsname{#7}%
    \addcontentsline{toc}{#1}{%
      \ifnum #2>\c@secnumdepth \else
        \protect\numberline{\csname the#1\endcsname.}%
      \fi
      #7}%
  \else
    \def\@svsechd{%
      #6{\hskip #3\relax
      \@svsec #8}%
      \csname #1mark\endcsname{#7}%
      \addcontentsline{toc}{#1}{%
        \ifnum #2>\c@secnumdepth \else
          \protect\numberline{\csname the#1\endcsname.}%
        \fi
        #7}}%
  \fi
  \@xsect{#5}}
\def\@seccntformat#1{\csname the#1\endcsname.\quad}
\makeatother

\newtheorem{theo}[equation]{Th\'eor\`eme}
\newtheorem{lem}[equation]{Lemme}

\newenvironment{remarque}{
\refstepcounter{equation}\trivlist%
\item[\hskip \labelsep{\bfseries Remarque \theequation.\ }]}%
{\endtrivlist}%

\makeatletter
\renewcommand\theequation{\thesection.\arabic{equation}}
\@addtoreset{equation}{section}
\makeatother
\newcommand{\carrenoir}{\rule{0.5em}{0.5em}}
\makeatletter
\newenvironment{demo}[1][\@empty]{\textbf{D\'emonstration~%
\ifx\@empty#1:\else #1~:\fi~}}
{\hfill\carrenoir\nolinebreak\vspace{2mm}}
\makeatother
\newcommand{\oper}[2]{\newcommand{#1}{\mathop{\mathrm{#2}}\nolimits} }

\oper{\scal}{scal}
\newcommand{\R}{\mathbb R}

\DeclareSymbolFont{greek}{OML}{ptmcm}{m}{it}
\DeclareMathSymbol{\codiff}{\mathord}{greek}{"0E}
\newcommand{\de}{\mathrm{ d }}

\title{Un théorème de la masse positive pour le problème de Yamabe en
dimension paire}
\author{Pierre Jammes}
\date{}
\begin{document}
\maketitle
{\small 
\textsc{Résumé.---}
Soit $(M,g)$ une variété compacte conformément plate de dimension $n\geq4$
et de courbure scalaire strictement positive. Selon une théorème de
la masse positive dû à Schoen et Yau, le terme constant dans le développement
de la fonction de Green du laplacien conforme est strictement positif
quand $(M,g)$ n'est pas conforme à la sphère ronde. Sur les variétés spin,
Ammann et Humbert en ont donné une démonstration élémentaire, basée 
sur une preuve de Witten. En utilisant les formes différentielles au lieu
des spineurs, nous en donnons une démonstration élémentaire sur
les variétés de dimension paire, sans autre hypothèse sur la topologie.

Mots-clefs : théorème de la masse positive, problème de Yamabe, formes
différentielles.

\medskip
\textsc{Abstract.---}
Let $(M,g)$ be a compact conformally flat manifold of dimension $n\geq4$
with positive scalar curvature. According to a positive mass theorem by 
Schoen and Yau, the constant term in the development of the Green function
of the conformal Laplacian is positive if $(M,g)$ is not conformally
equivalent to the sphere. On spin manifolds, there is an elementary proof 
of this fact by Ammann and Humbert, based on a proof of Witten. Using
differential forms instead of spinors, we give an elementary proof on
even dimensional manifolds, without any other topological assumption.


}

\section{Introduction}

Le problème de Yamabe consiste à trouver dans chaque classe conforme
d'une variété compacte une métrique dont la courbure scalaire est
constante. Après avoir décelé une erreur dans la solution \cite{ya60}
proposée par Yamabe, N.-S.~Trudinger a pu la corriger dans \cite{tr68}
pour les classes conformes n'admettant pas de métrique à courbure 
scalaire strictement positive. Les cas restants ont été résolu par
T.~Aubin \cite{au76} pour les variétés non conformément plate de
dimension $n\geq6$ et par R.~Schoen \cite{sc84} en dimension~3 à~5 et pour 
les variétés conformément plates (voir \cite{lp87}).

La démonstration de R.~Schoen consiste à se ramener à l'étude
d'une variété asymptotiquement euclidienne de la manière suivante:
si le laplacien conforme $L_g=\frac{4(n-1)}{n-2}\Delta_g+\scal_g$ est 
inversible (en particulier s'il existe une metrique de courbure scalaire
strictement positive dans la classe conforme de $g$), il existe pour
tout point $P\in M$ une fonction de Green $\Gamma_P$ telle que
$L_g\Gamma_P=\delta_P$ au sens des distributions. Comme cette fonction
admet un développement de la forme
\begin{equation}\label{intro:green}
\Gamma_P(x)=\frac1{4(n-1)\omega_{n-1}r^{n-2}}+A_P+f(x) 
\end{equation}
où $r=d(P,x)$, $f(x)=O(r)$ et $\omega_{n-1}$ désigne le volume de la sphère 
canonique $S^{n-1}$, la variété $(M\backslash\{P\},\Gamma_P^{\frac4{n-2}}g)$ 
est asymptotiquement euclidienne (cf.~\cite{lp87}). On appelle souvent 
\emph{projection stéréographique} cette construction car elle généralise 
la projection classique de $S^n\backslash\{P\}$ sur $\R^n$. 
Dans \cite{sc84}, R.~Schoen montre d'une part qu'on peut
résoudre le problème de Yamabe si $A>0$, et d'autre part que la constante
$A$ peut s'identifier à la masse de $(M\backslash\{P\},
\Gamma_P^{\frac4{n-2}}g)$, la masse étant un invariant riemannien des variétés 
asymptotiquement plates dont l'étude était initialement motivée par des
problèmes issus de la relativité générale.

On est donc ramené au problème de la masse positive, qui
consiste à montrer que la masse d'une variété asymptotiquement plate de
courbure scalaire positive est positive, l'annulation de la masse 
caractérisant $\R^n$. Il existe différents théorèmes de ce type dans 
des contextes plus ou moins généraux (voir le survol~\cite{he98}). Cette
conjecture a d'abord été montrée en dimension $3\leq n\leq7$ par Schoen et 
Yau dans \cite{sy79} et \cite{sy81}, et le cas des projections 
stéréographiques de variétés conformément plates est couvert par 
\cite{sy88}, ce qui achève la résolution du problème de Yamabe.

Par ailleurs E.~Witten a proposé dans \cite{wi81} une démonstration
du théorème de la masse positive sur les variétés asymptotiquement 
euclidiennes utilisant les spineurs (voir aussi \cite{ba86}). B.~Ammann et 
E.~Humbert l'ont simplifiée
dans le cas des projections stéréographiques ; sur les variétés
conformément plates leur démonstration est particulièrement concise et
élégante. Cependant, ces résultats nécessitent
l'hypothèse topologique que la variété soit spin. Le but de cet article
est de modifier la démonstration de \cite{ah05} de manière à s'affranchir
de cette condition. L'argument ne fonctionne qu'en dimension paire,
mais il n'y a aucune autre restriction topologique sur la variété,
pas même l'orientabilité:
\begin{theo}\label{intro:th}
Soit $(M,g)$ une variété riemannienne compacte et conformément plate,
de courbure scalaire strictement positive et de dimension $n\geq4$ paire. 
Pour tout $P\in M$, la masse $A_P$ est positive, et $A_P=0$ si et seulement 
si $(M,g)$ est conforme à la sphère canonique.
\end{theo}
La démonstration du théorème~\ref{intro:th} est assez similaire à celle
de \cite{ah05}. La principale différence consiste à raisonner sur des
formes différentielles au lieu des spineurs, et on utilisera le fait
que si $n$ es pair, la formule de Weitzenböck pour les $\frac n2$-formes
est assez simple et bien adaptée au problème.

Je remercie Emmanuel Humbert de m'avoir présenté la démonstration
de \cite{ah05}, ainsi que Vincent Minerbe de m'avoir recommandé la
lecture de \cite{lp87}.

\section{Démonstration}

 Commençons par quelques rappels techniques concernant le laplacien
agissant que les formes différentielles : il est défini par
$\Delta=\de\codiff+\codiff\de$ où la codifférentielle $\codiff$ agissant
sur les $p$-formes est l'adjoint $L^2$ de la différentielle $\de$.
Contrairement à une idée répandue, les opérateurs $\codiff$ et $\Delta$
sont bien défini quand $M$ n'est pas orientable.
Dans la suite, on dira qu'une forme différentielle $\omega$ est harmonique si
$\de\omega=\codiff\omega=0$ (ce qui est un condition plus forte que
$\Delta\omega=0$).
On peut aussi écrire localement $\codiff=(-1)^{n(p+1)+1}*\de*$
(le signe de la dualité de Hodge $*:\Omega^p(M)\to\Omega^{n-p}(M)$ dépend
d'un choix local d'orientation, mais cela n'affecte pas le signe de
$\codiff$), ce qui permet d'appliquer $\codiff$ à des formes
différentielles qui ne sont pas $L^2$. En outre, si la dimension $n$
est paire, l'opérateur $*:\Omega^{\frac n2}(M)\to\Omega^{\frac n2}(M)$
est conformément invariant.  En particulier, si
$\omega\in\Omega^{\frac n2}(M)$ alors
l'harmonicité de $\omega$ est une propriété conforme (contrairement à
$\Delta\omega=0$).

Considérons maintenant une variété compacte $(M^n,g)$ de dimension $n$ 
paire conformément plate et de courbure scalaire positive, et un point
$P\in M$. Comme la positivité de $A_P$ ne dépend pas du choix de la 
métrique dans la classe conforme de $g$ (\cite{au98}, proposition~5.41), on 
peut en fait supposer que la métrique $g$ est plate dans une petite
boule $B_P(r_0)$ de rayon $r_0$ et centrée en $P$. La première étape de 
la démonstration consiste à construire une $\frac n2$-forme harmonique 
sur $M\backslash\{P\}$, pour ensuite l'utiliser pour estimer la masse 
de la projection stéréographique de $M$. Cette étape nécessite des outils
spécifiques aux formes différentielles.

 On commence
par construire un modèle euclidien de forme harmonique : soit $\varphi_0\in
\Lambda^{\frac n2}\R^n$ une forme alternée qu'on identifie à la 
${\frac n2}$-forme différentielle invariante correspondante sur $\R^n$.
Cette forme est parallèle, donc $\de\varphi_0=\codiff\varphi_0=0$. 
L'inversion $i$ par rapport à la sphère unité est un difféomorphisme
conforme de $\R^n\backslash\{0\}$, on a donc aussi $\de i^*(\varphi_0)
=\codiff i^*(\varphi_0)=0$. Si on note $r$ la coordonnée radiale sur $\R^n$,
on a $i^*(g_\textrm{eucl})=r^{-4}g_\textrm{eucl}$, donc $|i^*(\varphi_0)|=
r^{-n}$ et de plus $\nabla_r(i^*(r^n\varphi_0))=0$ car la forme
$i^*(r^n\varphi_0)$ est constante le long d'une droite passant par l'origine.

\begin{lem}\label{demo:lem}
Soit $\varphi_0\in\Lambda^{\frac n2}\mathbb{R}^n$. Il existe une
${\frac n2}$-forme $\varphi$ harmonique sur $M\backslash\{P\}$
et telle que $\varphi=i^*(\varphi_0)+\varphi'$ sur $B_P(r_0)\backslash P$,
$\varphi'$ étant une forme lisse sur $B_P(r_0)$.
\end{lem}

\begin{demo}
On commence par choisir une fonction \emph{cut off} $\eta$ qui vaut 
1 sur $B_P(r_0/2)$ et 0 en dehors de $B_P(r_0)$, et on note $\bar\varphi$
la forme définie sur $M\backslash\{P\}$ par $\bar\varphi=\eta i^*(\varphi_0)$ 
sur $B_P(r_0)$ et prolongée par 0 en dehors. Comme $\de\bar\varphi=0$ 
sur $B_P(r_0/2)\backslash\{P\}$, forme différentielle $\de\bar\varphi$ se 
prolonge de manière lisse en une forme fermée sur $M$ à support dans 
$B_P(r_0)$ ; le 
lemme de Poincaré à support compact (\cite{bt82}, corollaire~4.7.1)
assure alors qu'il existe donc une 
${\frac n2}$-forme $\theta$ à support dans $B_P(r_0)$
telle que $\de\theta=\de\bar\varphi$ sur $B_P(r_0)\backslash\{P\}$. 
En particulier, $\de\bar\varphi$ est la restriction à $M\backslash\{P\}$ 
d'une forme exacte sur $M$.

 Par théorie de Hodge, on peut donc classiquement trouver une forme lisse
$\theta_1$ telle que $\de\theta_1=\de\bar\varphi$ sur $M\backslash\{P\}$
et $\codiff\theta_1=0$ (on utilise ici seulement le fait que $\codiff$ est
l'adjoint de $\de$, l'hypothèse d'orientation est superflue).
Comme on a aussi $\codiff\bar\varphi=0$ sur $B_P(r_0/2)\backslash\{P\}$,
le même argument permet de
trouver une forme $\theta_2$ lisse et fermée sur $M$ vérifiant
$\codiff\theta_2=\codiff\bar\varphi$. Il suffit alors de choisir
$\varphi=\bar\varphi-\theta_1-\theta_2$.
\end{demo}

\begin{remarque}
Dans \cite{ah05}, la construction d'un spineur harmonique 
sur $M\backslash\{P\}$ nécessite une condition sur la métrique, à savoir
que l'opérateur de Dirac est inversible (par exemple si la courbure scalaire
est strictement positive). Ce n'est pas le cas avec les formes 
différentielles : c'est l'utilisation du lemme de Poincaré qui permet de 
se dispenser de ce genre d'hypothèse.
\end{remarque}

\begin{demo}[du théorème \ref{intro:th}]
Si $\Gamma_P$ est la fonction de Green de $L_g$ au point $P$, on pose 
$G=4(n-1)\omega_{n-1}\Gamma_P$ et on pose $\tilde g=G^{\frac4{n-2}}g$.
Du fait que $L_g(G)=0$, la formule de transformation de la courbure scalaire 
pour les déformations conforme indique que $\scal_{\tilde g}=0$ sur 
$(M\backslash\{P\},\tilde g)$.

Comme la forme $\varphi$ fournie par le lemme~\ref{demo:lem} est de 
degré $\frac n2$, elle reste harmonique sur $(M\backslash\{P\},\tilde g)$,
et donc $\Delta\varphi=0$. 
J.-P.~Bourguignon a montré (\cite{bo81}, proposition~8.6 et remarque~8.7)
que pour les formes de degré~$\frac n2$ en dimension paire, le terme
de courbure de la formule de Weitzenböck s'exprime uniquement à l'aide
du tenseur de Weyl et de la courbure scalaire (voir aussi~\cite{la06},
proposition~4.2). Comme $M$ est conformément
plate et que $\scal_{\tilde g}=0$, on a simplement $\Delta\varphi=\nabla^*
\nabla\varphi=0$ sur $(M\backslash\{P\},\tilde g)$. Une intégration
par partie donne alors 
\begin{eqnarray}
\int_{M\backslash B_P(r)}|\nabla\varphi|_{\tilde g}^2\de v_{\tilde g}&=&
\int_{M\backslash B_P(r)}\langle\nabla^*\nabla\varphi,\varphi
\rangle_{\tilde g}\de v_{\tilde g}+\int_{S_P(r)}\langle\nabla_\nu\varphi,
\varphi\rangle_{\tilde g}\de s_{\tilde g}\nonumber\\
&=&\int_{S_P(r)}\langle\nabla_\nu\varphi,\varphi\rangle_{\tilde g}
\de s_{\tilde g}=\frac12\int_{S_P(r)}\partial_\nu|\varphi|_{\tilde g}^2
\de s_{\tilde g}\label{demo:eq1}
\end{eqnarray}
où $\nu$ est le vecteur unitaire (pour $\tilde g$) normal entrant de sphère
$S_P(r)$ et $\de s_{\tilde g}$ la forme volume induite sur
$S_P(r)$ par $\tilde g$.

D'une part, en utilisant le fait que sur les variétés conformément plates
la fonction $f$ dans (\ref{intro:green}) peut s'écrire $r\bar f(x)$ où
$\bar f(x)=O(1)$ quand $r\to0$  (\cite{lp87}, lemme~6.4), on a
\begin{eqnarray}
|\varphi|_{\tilde g}^2&=&G^{\frac{-2n}{n-2}}|\varphi|_g^2=\left(
\frac1{r^{n-2}}+4(n-1)\omega_{n-1}A+r\bar f(x)\right)^{\frac{-2n}{n-2}}
|i^*(\varphi_0)+\varphi'|_g^2\nonumber\\
&=&(1+4(n-1)\omega_{n-1}Ar^{n-2}+r^{n-1}\bar f(x))^{\frac{-2n}{n-2}}
|i^*(r^n\varphi_0)+r^n\varphi'|_g^2\nonumber\\
&=&(1+4(n-1)\omega_{n-1}Ar^{n-2}+r^{n-1}\bar f(x))^{\frac{-2n}{n-2}}\times
\nonumber\\
&&(1+r^n\langle i^*(r^n\varphi_0),\varphi'\rangle_g+r^{2n}|\varphi'|_g^2).
\end{eqnarray}
En utilisant le fait que $\nabla_r(i^*(r^n\varphi_0))=0$, on en déduit
quand $r\to0$ que
\begin{equation}\label{demo:eq2}
\frac\partial{\partial r}|\varphi|_{\tilde g}^2=-8n(n-1)\omega_{n-1}Ar^{n-3}
+o(r^{n-3})
\end{equation}

D'autre part, pour $r$ petit, on a aussi $\nu=-G^{-\frac2{n-2}}
\frac\partial{\partial r}\sim-r^2\frac\partial{\partial r}$
et $\de s_{\tilde g}=G^{\frac{2(n-1)}{n-2}}r^{n-1}\de s\sim
r^{-(n-1)}\de s$ où $\de s$ désigne la forme volume canonique de la sphère 
$S^{n-1}$. En conjonction avec (\ref{demo:eq1}) et (\ref{demo:eq2}), on
obtient
\begin{equation}
0\leq\int_{M\backslash B_P(r)}|\nabla\varphi|_{\tilde g}^2\de v_{\tilde g}
=\frac12\int_{S_P(r)}\partial_\nu|\varphi|_{\tilde g}^2\de s_{\tilde g}=
4n(n-1)\omega_{n-1}^2A+o(1)
\end{equation}
et on en déduit que $A\geq0$.

Il reste à traiter le cas où $A=0$. L'inégalité précédente permet d'affirmer
que pour toute forme $\varphi_0\in\Lambda^{\frac n2}\R^n$, la forme 
$\varphi$ correspondante sur $(M\backslash\{P\},\tilde g)$ est parallèle.
On a donc une base de $\frac n2$-formes différentielles parallèles sur
$(M\backslash\{P\},\tilde g)$, et on peut en déduire qu'il existe aussi
une base de champs de vecteurs parallèles : si on construit un champ 
de vecteur comme intersection de noyaux de formes parallèles, il sera
parallèle. Par conséquent la variété $(M\backslash\{P\},\tilde g)$ est 
plate, donc isométrique à $\R^n$, et $(M,g)$ est conforme à la sphère 
canonique.
\end{demo}

\noindent Pierre \textsc{Jammes}\\
Université d'Avignon et des pays de Vaucluse\\
Laboratoire d'analyse non linéaire et géométrie (EA 2151)\\
F-84018 Avignon\\
\texttt{Pierre.Jammes@univ-avignon.fr}

\end{document}